\documentclass{article}
\usepackage[english]{babel}
\usepackage[cp1250]{inputenc}
\usepackage{amsmath}
\usepackage{amsfonts}
\usepackage{amsthm}
\usepackage[mathcal]{euscript}
\usepackage{graphicx}
\usepackage{enumerate}
\newcommand{\cL}{{\cal L}}
\newcommand{\cC}{{\cal C}} 
\newcommand{\cD}{{\cal D}} 
\newcommand{\bbN}{\mathbb{N}}
\newcommand{\bbZ}{\mathbb{Z}}
\newcommand{\bbR}{\mathbb{R}}
\newcommand{\cB}{\Delta} 
\newcommand{\B}{\mathfrak{B}}
\newcommand{\W}{\mathfrak{W}}
\newtheorem{thm}{Theorem}[section]
\newtheorem{cla}[thm]{Claim}
\newtheorem{lem}[thm]{Lemma}
\newtheorem{prop}[thm]{Proposition}                   
\newtheorem{con}[thm]{Conjecture}
\newtheorem{cor}[thm]{Corollary}
\theoremstyle{definition}
\newtheorem{defi}[thm]{Definition}      
\def\inst#1{$^{#1}$}

\begin{document} 
 
\date{}
\title{Monochromatic triangles in two-colored plane}
\author{V\'{\i}t Jel\'{\i}nek\inst1 \and Jan Kyn\v cl\inst2 \and Rudolf Stola\v r\inst2 \and Tom\'a\v s Valla\inst2}

\maketitle

\begin{center}
{\footnotesize
\inst{1} Charles University, Faculty of Mathematics and Physics,\\
Institute for Theoretical Computer Science (ITI)\footnote{ITI is supported by project 1M0021620808 of the Czech Ministry of Education}\\
Malostransk\'e n\'am.~2/25, 118~00, Prague, Czech Republic \\
\texttt{jelinek@kam.mff.cuni.cz}
\\\ \\
\inst{2} Charles University, Faculty of Mathematics and Physics,\\
Department of Applied Mathematics (KAM)\footnote{KAM is supported by project MSM0021620838 of the Czech Ministry of Education}\\
Malostransk\'e n\'am.~2/25, 118~00, Prague, Czech Republic \\
\texttt{jankyncl@centrum.cz, riv@email.cz, valla@kam.mff.cuni.cz}
} \end{center}                                                                                

\begin{abstract}
We prove that for any partition of the plane into a closed set $C$ and an open 
set $O$ and for any configuration $T$ of three points, there is a  
translated and rotated copy of $T$ contained in $C$ or in $O$. 

Apart from that, we consider partitions of the plane into two sets whose 
common boundary is a union of piecewise linear curves. We show that for any such 
partition and any configuration $T$ which is a vertex set of a 
non-equilateral triangle there is a copy of $T$ contained in the 
interior of one of the two partition classes. Furthermore, we give the 
characterization of these ``polygonal'' partitions that avoid 
copies of a given equilateral triple.

These results support a conjecture of Erd\H os, Graham, Montgomery, 
Rothschild, Spencer and Straus, which states that every two-coloring of the 
plane contains a monochromatic copy of any nonequilateral triple of points; 
on the other hand, we disprove a stronger conjecture by the same authors, by 
providing non-trivial examples of two-colorings that avoid a given 
equilateral triple. 
\end{abstract}

\section{Introduction}          
Euclidean Ramsey theory addresses the problems of the following kind: 
assume that a finite configuration $X$ of points is given; for what values of 
$c$ and $d$ is it true that every coloring of the $d$-dimensional Euclidean 
space by $c$ colors contains a monochromatic congruent copy of $X$? The first 
systematic treatise on this theory appears in 1973 in a series of papers 
\cite{erdi,erdii,erdiii} by Erd\H os, Graham, Montgomery, Rothschild, Spencer 
and Straus. Since that time, many strong results have been obtained in this 
field, often related to high-dimensional configurations (see, e.g., 
\cite{fra,kri1,kri2,mat} or the survey \cite{gra}); however, there are basic 
`low-dimensional' problems that remain open.

In this paper, we consider the special case when $d=2$, $c=2$ and $|X|=3$; in 
other words, we study the configurations of three points in the Euclidean 
plane colored by two colors. We use the term \emph{triangle} to refer to any 
set of three points, including collinear triples of points, which we call 
\emph{degenerate} triangles. An $(a,b,c)$-triangle is a triangle whose edges, 
in anti-clockwise order, have respective lengths $a$, $b$ and~$c$. A 
$(1,1,1)$-triangle is also called a \emph{unit triangle}.

We say that a set of points $X\subseteq \bbR^2$ is a \emph{copy} of a set of 
points $Y\subseteq\bbR^2$, if $X$ can be obtained from $Y$ by translations and 
rotations in the plane. A \emph{coloring} is a partition of $\bbR^2$ into two 
sets $\B$ and $\W$. The elements of $\B$ and $\W$ are called \emph{black 
points} and \emph{white points}, respectively. We use the term \emph{boundary 
of $\chi$} to refer to the common boundary of the sets $\B$ and $\W$. Given a coloring 
$\chi=(\B,\W)$, we say that a set of points $X$ is \emph{monochromatic}, if 
$X\subseteq\B$ or $X\subseteq\W$. 

We say that a coloring $\chi$ \emph{contains} a triangle $T$, if there exists a 
monochromatic set $T'$ which is a copy of $T$; otherwise, we say that $\chi$ 
\emph{avoids} $T$. 

A coloring that avoids the unit triangle is easy to obtain: consider a 
coloring $\chi^*$ that partitions the plane into alternating half-open strips 
of width $\frac{\sqrt{3}}{2}$; formally, a point $(x,y)$ is black if an only 
if $n\sqrt{3}<y\le\left(n+\frac{1}{2}\right)\sqrt{3}$ for some integer $n$. 
It can be easily checked that $\chi^*$ avoids the unit triangle. We can even 
change the color of some of the points on the boundaries of the strips 
without creating any monochromatic unit triangle. Erd\H os et 
al.~\cite[Conjecture~1]{erdiii} have conjectured that this is essentially 
the only example of colorings avoiding a given triangle:

\begin{con}[Erd\H os et al.\ \cite{erdiii}]\label{con-silna}
For every triangle $T$ and every coloring $\chi$, if $\chi$ avoids 
 $T$, then $T$ is an equilateral $(l,l,l)$-triangle and $\chi$ is 
equal to an $l$-times scaled copy of the coloring $\chi^*$ defined above, up to  
possible modifications of the colors of the points on the boundary of the strips. 
\end{con}                                                                         

In Section~\ref{sec-poly} of this paper, we present a counterexample to this 
conjecture, and define a general class of colorings (which includes 
$\chi^*$ as a special case) that avoid the unit triangle.

On the other hand, the following conjecture by Erd\H os et al. \cite[Conjecture~3]{erdiii} 
remains open:

\begin{con}
[Erd\H os et al.\ \cite{erdiii}]\label{con-slaba} Every coloring $\chi$ 
contains every nonequilateral triangle $T$.
\end{con}                                                          
                                                                   
In the past, it has been shown that Conjecture~\ref{con-slaba} holds for 
special types of triangles $T$ (see, e.g., \cite{erdiii,gra,sha}). Our 
approach is different: we prove that the conjecture is valid for a restricted 
class of colorings $\chi$ and arbitrary $T$. In Section~\ref{sec-uzav}, we 
show that every coloring that partitions $\bbR^2$ into a closed set and an 
open set contains every triangle $T$. Then, in Section~\ref{sec-poly}, we 
consider \emph{polygonal} colorings, whose boundary is a union of piecewise 
linear curves (see page~\pageref{def-poly} for the precise definition). We 
show that Conjecture~\ref{con-slaba} holds for the polygonal colorings, but 
there are polygonal counterexamples to the stronger 
Conjecture~\ref{con-silna}. In fact, we are able to characterize all these 
polygonal counterexamples.

The following lemma from \cite{erdiii} offers a useful insight into the topic 
of monochromatic triangles in two-colored plane:

\begin{lem}
\label{lem-osm} Let $\chi$ be a coloring of the plane. The following holds:
\begin{enumerate}[(i)]
\item 
If $\chi$ contains an $(a,a,a)$-triangle for some $a>0$, then 
$\chi$ contains an $(a,b,c)$-triangle, for every $b,c>0$ such that 
$a,b,c$ satisfy the (possibly degenerate) triangle inequality.
\item
If $\chi$ contains an $(a,b,c)$-triangle, then $\chi$ contains an 
$(x,x,x)$-triangle for some $x\in\{a,b,c\}$.
\end{enumerate}
\end{lem}                                        
\begin{figure}
\begin{center}
\includegraphics[scale=0.9]{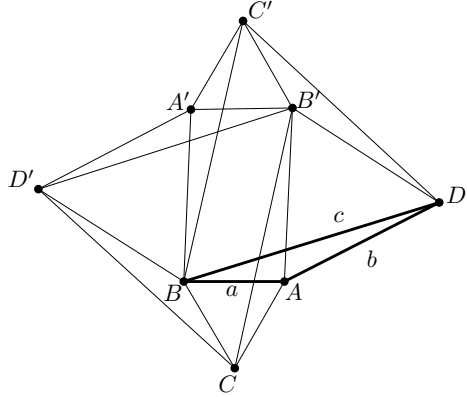}
\end{center} 
\caption[Proof of Lemma~\ref{lem-osm}]{The illustration of the proof of Lemma~\ref{lem-osm}}
\end{figure}\label{fig-osm}
\begin{proof}
The essence of the proof is the configuration in Figure~\ref{fig-osm}. The 
configuration consists of two $(a,a,a)$-triangles $ABC$ and $A'B'C'$, two 
$(b,b,b)$-triangles $ADB'$ and $A'D'B$ and two $(c,c,c)$-triangles $BDC'$ and 
$B'D'C$. To prove the first part of the lemma, assume, for a given $\chi$, 
that there is a monochromatic $(a,a,a)$-triangle $ABC$, and choose arbitrary 
$b$ and $c$ satisfying triangle inequality with $a$. Assume that $A$, $B$ and 
$C$ are all black. Furthermore, assume for contradiction that no 
$(a,b,c)$-triangle is monochromatic. Considering the configuration in 
Fig.~\ref{fig-osm}, we deduce that the points $B'$, $D$ and $D'$ are all 
white, otherwise one of the $(a,b,c)$-triangles $BAD$, $CAB'$ and $CBD'$ 
would be monochromatic. Then, $A'$ is black, due to $B'A'D'$, and $C'$ is 
white, due to $C'A'B$. It follows that $C'B'D$ is monochromatic, a 
contradiction.

The second part is proved by an analogous argument: assume that $BAD$ is an
all-white monochromatic triangle and that the statement does not hold. Then $B'$,
$C$ and $C'$ are all black, due to $ADB'$, $ABC$ and $BDC'$. $A'$ is white, due
to $A'B'C'$; $D'$ is black, due to $A'D'B$, and $B'D'C$ is monochromatic. 

This concludes the proof.
\end{proof}                                                          

From Lemma~\ref{lem-osm}, we obtain directly the following facts:

\begin{cor}\label{cor-osm} For every coloring $\chi$ the following holds:

\begin{enumerate}[(i)]
\item $\chi$ contains every triangle if and only if $\chi$
contains every equilateral triangle.
\item $\chi$ contains every non-equilateral triangle if
and only if there is an $a_0>0$ such that $\chi$ contains the equilateral
$(a,a,a)$-triangle for all values of $a>0$ different from $a_0$.
\item $\chi$ contains an $(a,b,c)$-triangle if and only if $\chi$
contains a $(b,a,c)$-triangle.
\end{enumerate}
\end{cor}

\section{Coloring by closed and open sets}\label{sec-uzav}  

The aim of this section is to prove the following result:

\begin{thm}\label{thm-uzav} 
Let $\chi=(\B,\W)$ be a coloring such that $\B$ is closed and $\W$ is open. 
Then $\chi$ contains every triangle $T$.
\end{thm}

By Corollary~\ref{cor-osm}, it suffices to prove Theorem~\ref{thm-uzav} for 
the case when $T$ is an arbitrary equilateral triangle. Moreover, since scaling 
does not affect the topological properties of $\B$ and $\W$, we only need to 
consider the case when $T$ is the unit triangle. Before stating the proof, we 
introduce a definition and prove an auxiliary result.

\begin{defi} Let $\varepsilon>0$. An $(a,b,c)$-triangle whose edge-lengths
satisfy $1-\varepsilon\le a,b,c \le 1+\varepsilon$ is called \emph{an
$\varepsilon$-almost unit triangle}.
\end{defi}

Suppose that an orthogonal coordinate system is given in the plane. For 
$a>0$, let $Q(a)$ be the closed square with vertices 
${(a,a),(-a,a),(-a,-a),(a,-a)}$. 

\begin{prop}\label{prop-almost} 
Let $Q(3)=\B \cup \W$ be a decomposition of the square $Q(3)$ into two 
disjoint sets such that there is no monochromatic unit triangle in $Q(3)$. 
Then for every $\varepsilon > 0$ both $\B$ and $\W$ contain an 
$\varepsilon$-almost unit triangle.
\end{prop}

\begin{figure}
\begin{center}\includegraphics{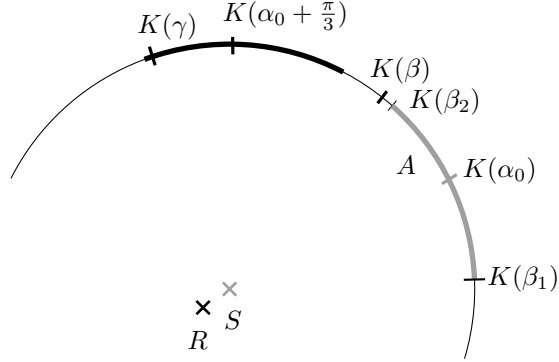}\end{center}
\caption{Illustration of the proof of Proposition~\ref{prop-almost}}\label{fig-uzav}  
\end{figure}

\begin{proof}
Let $\varepsilon$ be a given positive number. Assume that 
we are given a partition $\B\cup\W=Q(3)$ such that $Q(3)$ does not contain any 
monochromatic unit triangle. For contradiction, assume that one of the classes, 
wlog the class $\B$, does not contain any $\varepsilon$-almost unit triangle.

There is a white point $S$ and a black point $R$ in $Q(1)$ such that $|R-S| < 
\varepsilon$ (otherwise the whole square $Q(1)$ would be monochromatic). Let 
$\cC$ be the unit circle centered at $S$. For every $\alpha \in \bbR$, let 
$K(\alpha)$ denote the point of $\cC$ with coordinates $(x_S+\cos(\alpha), 
y_S+\sin(\alpha))$, where $(x_S,y_S)$ are the coordinates of $S$.

Note that the distance between $R$ and any point on $\cC$ is always in the 
interval $(1-\varepsilon, 1+\varepsilon)$; thus, for every $\alpha$, the 
points $K(\alpha)$ and $K(\alpha + \frac{\pi}{3})$ must have different 
colors, otherwise they would form a monochromatic white unit triangle with 
$S$ or a monochromatic black $\varepsilon$-almost unit triangle with $R$.

Let $K(\alpha_0)$ be a white point, then $K(\alpha_0 + \frac{\pi}{3})$ is 
black (see Fig.~\ref{fig-uzav}). Note that for every $\alpha \in (\alpha_0-\varepsilon, 
\alpha_0+\varepsilon)$ the distance between $K(\alpha)$ and $K(\alpha_0 + 
\frac{\pi}{3})$ is in the interval $(1-\varepsilon,1+\varepsilon)$, so the 
whole arc $\{K(\alpha); \alpha \in (\alpha_0-\varepsilon, 
\alpha_0+\varepsilon)\}$ is white. Let $A=\{K(\alpha); \alpha \in (\beta_1, 
\beta_2)\}$ be the maximal open white arc of $\cC$ containing the point 
$K(\alpha_0)$. Then the whole arc $\{K(\alpha); \alpha \in 
(\beta_1+\frac{\pi}{3},\beta_2+\frac{\pi}{3})\}$ is black. By definition of 
$A$, there exists $\beta \in (\beta_2, \beta_2 + \frac{\varepsilon}{2})$ such 
that $K(\beta)$ is black. There also exists $\gamma \in (\beta_2+\frac{\pi}{3}
 - \frac{\varepsilon}{2}, \beta_2+\frac{\pi}{3})$ such that $K(\gamma)$ is 
black. But then $(\gamma-\beta) \in 
(\frac{\pi}{3}-\varepsilon,\frac{\pi}{3})$, so the distance between the black 
points $K(\beta)$ and $K(\gamma)$ is in the interval $(1-\varepsilon, 1)$, 
hence the three points $R, K(\beta),K(\gamma)$ form a black 
$\varepsilon$-almost unit triangle\thinspace ---\thinspace a contradiction.
\end{proof}                   

We are now ready to prove the main result of this section.

\begin{proof}[Proof of Theorem~\ref{thm-uzav}]            
Let $\chi=(\B,\W)$ be a coloring, with $\B$ closed. By 
Corollary~\ref{cor-osm}, it is sufficient to show that $\chi$ contains the 
unit triangle. Assume, for contradiction, that this is not the case. Let 
$\B_0 = Q(3) \cap \B$ and let $\W_0=Q(3)\cap\W$. Clearly, neither $\B_0$ nor 
$\W_0$ contain the unit triangle, so by Proposition~\ref{prop-almost}, both 
these sets contain $\varepsilon$-almost unit triangles for every 
$\varepsilon>0$. In particular, the set $\B_0$ contains, for every 
$n\in\bbN$, a $\frac{1}{n}$-almost unit triangle $X_nY_nZ_n$. 
 
Since $\B_0$ is a compact set, the set $\B_0^3 = \B_0\times\B_0\times\B_0$ is compact as well. The 
sequence $\{(X_n,Y_n,Z_n); n \in \bbN\}$ is an infinite sequence of points in 
$\B_0^3$, so there exists a convergent subsequence 
$\{(X_{n_k},Y_{n_k},Z_{n_k}); k \in \bbN\}$. Let $(X,Y,Z) \in \B_0^3$ be its 
limit. Then $X,Y,Z\in \B$ are limits of the sequences $\{X_{n_k};k \in \bbN\}$, 
$\{Y_{n_k};k \in \bbN\}$, and $\{Z_{n_k};k \in \bbN\}$, respectively. The 
Euclidean distance is a continuous function of two variables, so $|X-Y| = 
\lim_{k \to \infty} |X_{n_k}-Y_{n_k}|=1$, similarly $|Y-Z|=|Z-X|=1$. Thus, 
$\{X,Y,Z\}$ is a black unit triangle in $Q(3)$, which is a contradiction.
\end{proof}

\section{Polygonal colorings}\label{sec-poly}
Throughout this section, $\cC(A)$ denotes the unit circle with center 
$A$, and $\cD(A)$ denotes the closed unit disc with center $A$. 

In this section, we consider \emph{polygonal} colorings of the
plane, defined as follows:
\begin{defi}\label{def-poly} A coloring $\chi=(\B,\W)$ is said to be \emph{polygonal},
if it satisfies the following conditions (see an example in Fig.~\ref{fig-poly}):
\begin{figure}[ht]
\begin{center}\includegraphics{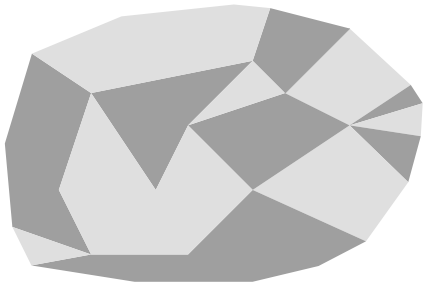}\end{center}
\caption{Example of a polygonal coloring}\label{fig-poly}
\end{figure}
\begin{itemize}
\item Each of the two sets $\B$ and $\W$ is contained in the
closure of its interior.
\item
The boundary of $\chi$ (denoted by $\cB$) is a union of straight line 
segments (called \emph{boundary segments}). Two boundary segments may only 
intersect at their endpoints. We allow these segments to be unbounded, i.e., 
a boundary segment may in fact be a half-line or a line. An endpoint of a 
boundary segment is called a \emph{boundary vertex}. We may assume that if 
exactly two boundary segments meet at a boundary vertex, then the two 
segments do not form a straight angle, because otherwise they
could be replaced with a single boundary segment. Note that with this condition, the 
boundary segments and boundary vertices of $\chi$ are determined uniquely.
\item
Every bounded region of the plane is intersected by only finitely many 
boundary segments (which implies that every bounded region contains only 
finitely many boundary vertices).
\end{itemize}
\end{defi}           
Note that these conditions imply that a sufficiently small disc around an interior
point of a boundary segment is separated by the boundary segment
into two halves, one of which is colored black and the other white. Note also 
that we make no assumptions about the colors of the points on the boundary~$\cB$.

We say that a coloring $\chi'$ is a \emph{twin} of a coloring $\chi$ if the 
two colorings have the same boundary and they assign the same colors to the 
points outside this boundary.

\begin{figure}[htb]
    \begin{center}
        \includegraphics[scale=0.7]{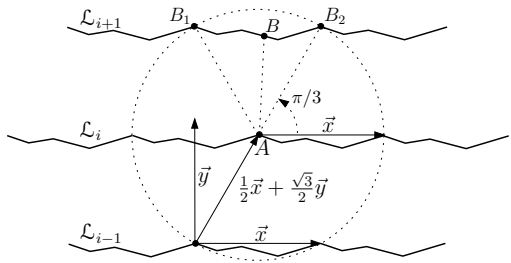}
    \end{center}
    \caption{The boundary of a zebra-like coloring}
    \label{fig:except}
\end{figure}

The main aim of this section is to prove that every polygonal coloring 
contains every nonequilateral triangle, and to characterize the polygonal 
colorings that avoid an equilateral triangle. To achieve this, we need the 
following definition:

\begin{defi}\label{def-strip} A coloring $\chi=(\B,\W)$ is called
\emph{zebra-like} if it has the following form:
the boundary of $\chi$ is a disjoint union of infinitely many
continuous curves $\cL_i; i\in\bbZ$ with the following
properties (see Fig.~\ref{fig:except}):
\begin{enumerate}[(a)]
\item There is a unit vector $\vec x$ such that for every $i\in\bbZ$, $\cL_i+\vec x=\cL_i$.
In other words, the $\cL_i$ are invariant upon a translation of length 1.
\item For every $i\in\bbZ$, the curve $\cL_{i+1}$ is a translated copy of $\cL_i$.
Moreover, there is a unit vector $\vec y$ orthogonal to $\vec x$, so that
\[
\cL_{i+1}=\cL_i+\frac{1}{2}\vec x+\frac{\sqrt{3}}{2}\vec y.
\]
In other words, for an arbitrary boundary point $X\in\cL_i$, the points $Y=X+\vec x$
and $Z=X+\frac{1}{2}\vec x+\frac{\sqrt{3}}{2}\vec y$ belong to the boundary as well.
Note that $XYZ$ is a unit triangle, and that $Y\in\cL_i$ and $Z\in\cL_{i+1}$. 
\item
For every $i\in\bbZ$, the interior of the region delimited by 
$\cL_i\cup\cL_{i+1}$ is colored with a different color than the interior of 
the region delimited by $\cL_{i-1}\cup\cL_{i}$.
\item
For two points $A$ and $B$, let $\theta_{AB}$ denote the size of
the acute angle formed by the segment $AB$ and the vector $\vec
x$. For every $i\in\bbZ$ and every two points $A\in\cL_i$ and
$B\in\cL_{i+1}$, the following holds: $\|AB\|>1$ if and only if
$\theta_{AB}<\frac{\pi}{3}$.

This last condition can also be stated in the following equivalent form:
Let $A\in\cL_i$ be an arbitrary point on the boundary. Let $B_1=A-\frac{1}{2}\vec
x+\frac{\sqrt{3}}{2}\vec y$ and $B_2=A+\frac{1}{2}\vec x+\frac{\sqrt{3}}{2}\vec y$
(the two points $B_1, B_2$ belong to $\cL_{i+1}$ by the previous conditions), and let
$A'=A+\sqrt{3}\vec y$ (so that $A'\in\cL_{i+2}$). Under these assumptions, the
portion of $\cL_{i+1}$ between $B_1$ and $B_2$ is contained inside of the closed lens-shaped
region $\cD(A)\cap\cD(A')$ and no other point of $\cL_{i+1}$ is inside this region.
\end{enumerate}
\end{defi}          

We stress that a zebra-like coloring is not necessarily polygonal.

\subsection{The result}
The following theorem is the main result of this section:

\begin{thm}\label{thm-poly}
For a polygonal coloring $\chi$, the following conditions are equivalent:
\begin{enumerate}
\item[(C1)] The coloring $\chi$ is a zebra-like polygonal coloring.  
\item[(C2)]
The coloring $\chi$ has a twin $\chi'$ which avoids the unit triangle.
\item[(C3)]
For every monochromatic unit triangle $ABC$, at least one of the three points
$A,B$ and $C$ belongs to the boundary of $\chi$.                               
\end{enumerate}
\end{thm}

Clearly, the condition (C2) of Theorem~\ref{thm-poly} implies the
condition (C3), so we only need to prove that (C1) implies (C2) and that
(C3) implies (C1).

The proof is organized as follows: we first prove that (C3)$\Rightarrow$(C1). 
This part of the proof proceeds in several steps: first of all, we use 
the condition (C3) to describe the set $\cB(\chi)\cap\cC(A)$, 
where $A$ is a boundary point. Then we apply a continuity argument 
to extend this information into a global description of~$\chi$.  

Next, in Theorem~\ref{thm-obarv}, we prove that every (not necessarily 
polygonal) zebra-like coloring has a twin that avoids the unit triangle, 
which shows that (C1)$\Rightarrow$(C2), completing the proof of 
Theorem~\ref{thm-poly}. 

In the last part of this section, we show  
that Theorem~\ref{thm-poly} implies that every polygonal coloring contains a 
monochromatic copy $T$ of a given non-equilateral triangle, with the vertices 
of $T$ avoiding the boundary. 

\subsection{The proof}
We begin with an auxiliary lemma:
\begin{lem}\label{lem-triprimky}
Let $q_1, q_2, q_3$ be (not necessarily distinct) lines in the
plane, not all three parallel. Then exactly one of the following
possibilities holds:
\begin{enumerate}
\item The lines $q_1, q_2, q_3$ intersect at a common point and every two of them form an angle
    $\frac{\pi}{3}$.
\item
There exist only finitely many unit triangles $ABC$ such that $A \in q_1$, $B
\in q_2$ and $C \in q_3$.
\end{enumerate}
\end{lem}

\begin{proof}
It can be easily checked that the two conditions cannot hold
simultaneously: in fact, if the three lines satisfy the first
condition, then for every point $A\in q_1$ whose distance from
the other two lines is at most 1 there are points $B\in q_2$ and
$C\in q_3$ such that $ABC$ is a unit triangle. We now show that at
least one of the two conditions holds.

Since the three lines are not all parallel, we may assume that neither $q_1$ 
nor $q_2$ is parallel to $q_3$. Consider a Cartesian coordinate system whose 
$y$\nobreakdash-axis is $q_3$. There exist real numbers $a_1, a_2, b_1, b_2$ 
such that for $i \in \{1,2\}$ we have $q_i=\{(x,y)\in \bbR^2; y=a_ix+b_i\}$. 
Let $ABC$ be a unit triangle with $A=(x_1,y_1) \in q_1$, $B=(x_2,y_2) \in q_2$ 
and $C\in q_3$, and assume that $A,B,C$ are in the counter-clockwise order 
(the other case is symmetric). Then $C=(\frac{x_1+x_2}{2},\frac{y_1+y_2} 
{2})+\frac{\sqrt{3}}{2}(y_1-y_2, x_2-x_1)$. The point $C$ lies on $q_3$, 
which implies the following equality:

\begin{equation}
    \frac{x_1+x_2}{2}+\frac{\sqrt{3}(y_1-y_2)}{2}=0 \label{equa1}
\end{equation}

Points $A$ and $B$ are at the distance $1$, from which we get

\begin{equation}
    (x_1-x_2)^2+(y_1-y_2)^2=1 \label{equa2}
\end{equation}

By combining \eqref{equa1} and \eqref{equa2} and eliminating $y_1, y_2$ we get

\[
\left(\frac{x_1+x_2}{2}\right)^2={\frac{3}{4}}\left(1-(x_1-x_2)^2\right),
\]
which yields
\begin{equation}
    x_1^2+x_2^2-x_1x_2={\frac{3}{4}}. \label{equa3}
\end{equation}

Substituting $y_1=a_1x_1+b_1$ and $y_2=a_2x_2+b_2$ into \eqref{equa1} gives

\begin{equation}
    {\frac{1+\sqrt{3}a_1}{2}}x_1 + {\frac{1-\sqrt{3}a_2}{2}}x_2 +
    {\frac{\sqrt{3}}{2}}(b_1-b_2)=0 \label{equa4}
\end{equation}
If both $\frac{1+\sqrt{3}a_1}{2}$ and $\frac{1-\sqrt{3}a_2}{2}$
are equal to zero, then the equality \eqref{equa4} degenerates and
we get that $a_1=-\frac{1}{\sqrt{3}}$, $a_2=\frac{1}{\sqrt{3}}$
and $b_1=b_2$, so the first case of the statement holds.

In the other case, suppose (wlog) that $\frac{1+\sqrt{3}a_1}{2}
\neq 0$. From \eqref{equa4} we can obtain that $x_1=cx_2+d$ for
some reals $c,d$. By substituting it into \eqref{equa3} we get a
quadratic equation for the variable $x_2$, where the leading
coefficient is equal to $c^2-c+1=(c-\frac{1}{2})^2+\frac{3}{4} >
0$, so there exist at most two possible values for $x_2$, thus at
most two possible locations of $B$ and at most four possible unit
triangles $ABC$.
\end{proof}

Throughout the rest of this section, we assume that $\chi$ is a fixed 
polygonal coloring satisfying the condition (C3) of Theorem~\ref{thm-poly}. 
Every boundary segment can be regarded as a common edge of two (possibly 
unbounded) polygonal regions, one of which is white and the other black. We 
choose an orientation of the boundary segments in the following way: a 
boundary segment with endpoints $A$ and $B$ is directed from $A$ to $B$ if 
the white region adjacent to this segment is on the left hand side from the 
point of view of an observer walking from $A$ to~$B$.

\begin{defi}\label{def-vhodny}
A boundary point $A\in\cB$ is called \emph{feasible}, if $A$ is not a boundary
vertex, and the unit circle $\cC(A)$ does not contain any boundary vertex.
An \emph{infeasible} point is a point on the boundary that is not feasible.
\end{defi}

We may easily see that every bounded subset of the plane contains only 
finitely many infeasible points.

The first step in the proof of the main result is the description of the set of 
all the boundary points at the unit distance from a given feasible point $A$.

Let $A$ be a fixed feasible point, let $s$ be the boundary segment containing 
$A$. The set $\cB\cap\cC(A)$ is finite, by the definition of polygonal 
coloring; on the other hand, this set is nonempty, otherwise we could find 
two points $B,C$ of $\cC(A)$ such that $ABC$ is a unit triangle, with $B$ and 
$C$ in the interior of the same color class. By shifting the triangle $ABC$ 
slightly in a suitable direction, we would obtain a monochromatic unit 
triangle avoiding the boundary, which is forbidden by the condition (C3).

In the following arguments, we will use a Cartesian coordinate system whose 
origin is the point $A$, and whose $x$-axis is parallel to $s$ and has the 
same orientation. We shall assume that the $x$-axis and the segment $s$ is 
directed left-to-right and the $y$-axis is directed bottom-to-top. Assuming 
this coordinate system, we let $P(\alpha,A)$ denote the point of $\cC(A)$ 
with coordinates $\left(\cos(\alpha),\sin(\alpha)\right)$. If no ambiguity 
arises, we write $P(\alpha)$ instead of $P(\alpha, A)$.

\begin{lem}
Let $B=P(\alpha)$ be an arbitrary element of $\cB\cap\cC(A)$, let $t$ be the boundary
segment containing $B$ (the segment $t$ is determined uniquely, because $A$ is
a feasible point).
Then the segments $s$ and $t$ are parallel.
\end{lem}
\begin{proof}
For contradiction, assume that $s$ and $t$ are not parallel, let
$\sigma\in(0,\pi)$ be the angular slope of $t$ with respect to the
coordinate system established above, i.e., $\sigma$ is the angle
formed by the lines containing $s$ and $t$.

First of all, note that the point $C=P(\alpha+\frac{\pi}{3})$ lies on the
boundary $\cB$; otherwise, a sufficiently small translation of the unit
triangle $ABC$ in a suitable direction  would yield a counterexample to
condition (C3) (here we use the assumption that $s$ and $t$ are not parallel).
Let $u$ be the boundary segment containing $C$, and let $\tau$ be the angular
slope of $u$.

Secondly, we may deduce that
$\{\sigma,\tau\}=\{\frac{\pi}{3},\frac{2\pi}{3}\}$, and the three
lines containing $s$, $t$ and $u$ all meet at one point. If this
were not the case, then by Lemma~\ref{lem-triprimky} there would
be only finitely many unit triangles with vertices belonging to
the three segments $s$, $t$ and $u$. Thus, we could find a unit
triangle $A'B'C'$ with $A'\in s$, $B'\in t$ and $C'\not\in\cB$,
which is impossible, by the argument presented in the previous
paragraph. By repeating this argument with
$\{\alpha+\frac{i\pi}{3};\ i=1,\dots,5\}$ in place of
$\alpha$, we obtain the following conclusions:
\begin{itemize}
\item 
The six points $\{P(\alpha+\frac{i\pi}{3});\ i=1,\dotsc,6\}$ all belong to 
the boundary $\cB$.
\item 
The lines passing through the boundary segments containing these six points 
all meet at one point.
\item 
The boundary segments containing $P(\alpha)$, $P(\alpha+\frac{2\pi}{3})$ and 
$P(\alpha+\frac{4\pi}{3})$ all have  the same slope.
\end{itemize}
This is a contradiction, because three parallel segments intersecting a circle
in three distinct points cannot belong to a single line, and two
parallel lines do not intersect.
\end{proof}

\begin{figure}[ht!] \begin{center} \includegraphics{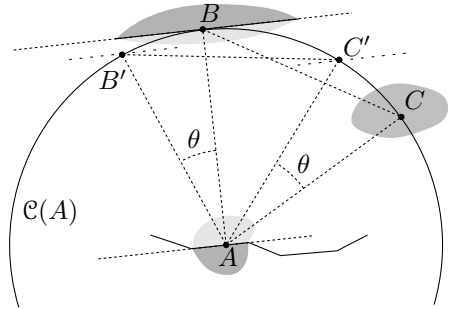} \end{center}
    \caption{Illustration of the proof of Lemma~\ref{lem-pipul}}
    \label{fig-pipul} \end{figure}

\begin{lem}\label{lem-pipul}
$P(\frac{\pi}{2})\not\in\cB$, $P(-\frac{\pi}{2})\not\in\cB$.
\end{lem}
\begin{proof}
For contradiction, assume that $B=P(\frac{\pi}{2})\in\cB$ (the case of 
$P(-\frac{\pi}{2})$ is symmetric), let $t$ denote the boundary segment 
containing $B$. Let $C=P(\frac{\pi}{6})$. We distinguish the following cases:
\begin{itemize}
\item
The segment $t$ has the same orientation as the segment $s$. In
this case, by applying a rotation around the center $C$ and then,
if $C\in\cB$, a suitable translation, we may transform the triple
$ABC$ into a monochromatic triple with vertices avoiding the
boundary, contradicting (C3).
\item
The segments $s$ and $t$ have opposite orientations (i.e., $t$ is
oriented right-to-left, which means that there is a white region
touching $t$ from below); furthermore, either $C\in\cB$ or $C$ is
in the interior of the white color. In such case, we may rotate
the configuration $ABC$ around the center of the segment $AB$ to
obtain a unit triangle in the interior of the white color.
\item
The segments $s$ and $t$ have opposite orientations and the point
$C$ is in the interior of the black color. Let $\theta$ be the
maximal angle with the properties that for every $\alpha\in
(\frac{\pi}{2},\frac{\pi}{2}+\theta)$ the point $P(\alpha)$ lies in the
interior of the white color and for every
$\alpha\in(\frac{\pi}{6},\frac{\pi}{6}+\theta)$ the point $P(\alpha)$ lies in
the interior of the black color. The value of $\theta$ is well defined, and by
the previous assumptions, $0<\theta<\frac{\pi}{3}$. Let
$B'=P(\frac{\pi}{2}+\theta)$ and $C'=P(\frac{\pi}{6}+\theta)$. By the maximality
of $\theta$, at least one of the two points lies on the boundary, and the
boundary segment passing through this point is directed left-to-right (see
Fig.~\ref{fig-pipul}). As in the first case of this proof, we may rotate and
translate the configuration $AB'C'$ to obtain a monochromatic unit triangle.
\end{itemize}

In all cases we get a contradiction.  
\end{proof}

The previous two lemmas imply that if $A$ is a feasible point,
then no boundary segment is tangent to $\cC(A)$.

\begin{lem}\label{lem-bilapod}
Let $B=P(\alpha)\in\cB$ be a point on the boundary, let $t$ be the
boundary segment containing this point. If $\alpha\in
(\frac{\pi}{6}, \frac{5\pi}{6})$ or
$\alpha\in(\frac{7\pi}{6},\frac{11\pi}{6})$, then $s$ and $t$ have
opposite orientation.
If $|\alpha|<\frac{\pi}{6}$ or $|\alpha-\pi|<\frac{\pi}{6}$,
then $s$ and $t$ have the same orientation.
\end{lem}

\begin{proof}
We first consider the case $\alpha\in(\frac{\pi}{6}, \frac{5\pi}{6})$
or $\alpha\in(\frac{7\pi}{6},\frac{11\pi}{6})$.
The proof is analogous to the proof of the first part of 
Lemma~\ref{lem-pipul}: if $t$ had the same orientation as $s$, we could take 
$C=P(\frac{\pi}{3}+\alpha)$ and then by rotating and translating the unit 
triangle $ABC$ we would get a contradiction. Note that the condition 
$\alpha\in(\frac{\pi}{6}, 
\frac{5\pi}{6})\cup(\frac{7\pi}{6},\frac{11\pi}{6})$ guarantees that $C$ is 
either the leftmost or the rightmost point of the triangle $ABC$, so whenever 
we start rotating the triangle $ABC$ around $C$, the two points $A,B$ move 
into the interior of the same color.

The case $|\alpha|<\frac{\pi}{6}$ or $|\alpha-\pi|<\frac{\pi}{6}$ can be proven analogously.
\end{proof}

\begin{figure}[Ht!]
\begin{center}
\includegraphics{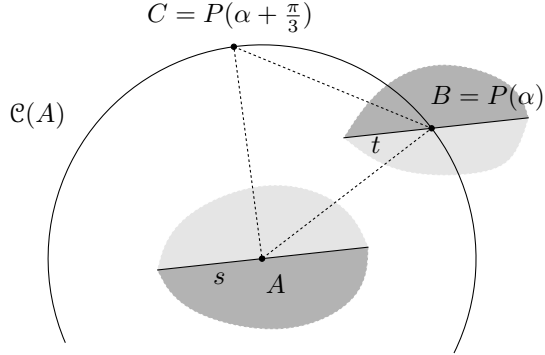}
\end{center}
\caption{Illustration of the proof of Lemma~\ref{lem-pitre}}
\label{fig-pitre}
\end{figure}

\begin{lem}\label{lem-pitre}
$P(\alpha)\in\cB$ if and only if $P(\alpha+\frac{\pi}{3})\in\cB$.
\end{lem}
\begin{proof}
It suffices to prove one implication, the other case is symmetric.
Assume that for some $\alpha$ we have $P(\alpha)\in\cB$ and
$P(\frac{\pi}{3}+\alpha)\not\in\cB$. Let $B=P(\alpha)$,
$C=P(\frac{\pi}{3}+\alpha)$, and let $t$ be the boundary segment
containing $B$. We consider the following cases:
\begin{itemize}
\item
If $s$ and $t$ have opposite orientation, we may rotate $ABC$
around the center of $AB$ to obtain a monochromatic unit triangle
in the interior of one color (see Fig.~\ref{fig-pitre}). Here we
use the fact that $\alpha\neq\frac{\pi}{2}$, which follows from
Lemma~\ref{lem-pipul}.
\item
If $s$ and $t$ have the same orientation, a small translation in a
suitable direction transforms $ABC$ into a monochromatic unit
triangle.
\end{itemize}
In both cases we get a contradiction.
\end{proof}

\begin{lem}
For every $\theta$ there is exactly one value of
$\alpha\in[\theta,\theta+\frac{\pi}{3})$ such that $P(\alpha)\in\cB$.
\end{lem}
\begin{proof}
By Lemma~\ref{lem-pitre}, if the statement holds for some value of $\theta$, it
holds for all other values of $\theta$ as well. Thus, it is enough to prove the
lemma for $\theta=\frac{\pi}{2}$.

Clearly, there is at least one $\alpha\in[\frac{\pi}{2},\frac{5\pi}{6})$ such 
that $P(\alpha)\in\cB$; otherwise, the set $\cC(A)\cap\cB$ would be empty, 
which is impossible.

Assume that there are $\alpha$ and $\alpha'$ such that
$\frac{\pi}{2}\le\alpha<\alpha'<\frac{5\pi}{6}$ with
$P(\alpha)\in\cB$ and $P(\alpha')\in\cB$. Let us fix $\alpha$ and
$\alpha'$ as small as possible. Let $t$ and $t'$ be the boundary
segments containing $P(\alpha)$ and $P(\alpha')$. The circle
$\cC(A)$ consists of alternating black and white arcs and one of
these arcs has $P(\alpha)$ and $P(\alpha')$ for endpoints. It
follows that one of the segments $t$, $t'$ has the same
orientation as the segment~$s$, contradicting
Lemma~\ref{lem-bilapod}.
\end{proof}

\begin{figure}[Ht!]
\begin{center}
\includegraphics[scale=.9]{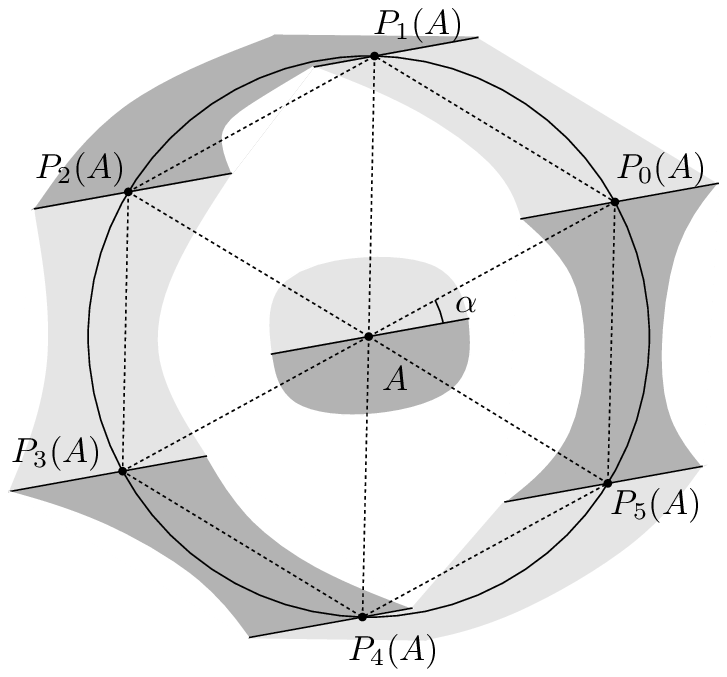}
\end{center}
\caption{Illustration of Claim~\ref{tvrz-sum}}
\label{fig-sum}
\end{figure}

Before we proceed with the proof of the main result, we summarize the lemmas proved
so far (and introduce some related notation) in the following
claim (see Fig.~\ref{fig-sum}):

\begin{cla}\label{tvrz-sum}
Let $A\in\cB$ be an arbitrary feasible point. The circle $\cC(A)$ intersects the
boundary $\cB$ at exactly six points, which form the vertex set of a
regular hexagon. These six points will be denoted by $P_0(A),\dotsc,P_5(A)$,
where $P_i(A)= P(\alpha+\frac{i\pi}{3},A)$ with $\alpha\in\left( -\frac{\pi}{6},\frac{\pi}{6} \right)$
(this determines $P_i(A)$ uniquely). The boundary segments containing the six
points $P_i(A)$ are all parallel to the boundary segment $s$ containing the point $A$.
The boundary segments containing the points  $P_0(A)$ and $P_3(A)$ have the
same orientation as $s$, whereas the boundary segments containing $P_1(A)$, $P_2(A)$,
$P_4(A)$ and $P_5(A)$ have opposite orientation.
\end{cla}

Now we use Claim~\ref{tvrz-sum} to get more global information about the boundary.

\begin{lem}
\label{lem-uhel}
Let $u_1$ and $u_2$ be two boundary segments that share a common endpoint $X$.
The size of the convex angle formed by these two segments is greater than $\frac{2\pi}{3}$.
\end{lem}                                                    

\begin{figure}[ht]
\begin{center}\includegraphics{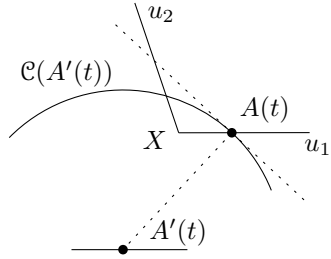}\end{center}
\caption{Illustration of the proof of Lemma~\ref{lem-uhel}}\label{fig-uhel}
\end{figure}

\begin{proof}
For contradiction, assume that for some $u_1$, $u_2$ and $X$, the
statement of the lemma does not hold (see Fig.~\ref{fig-uhel}). We may assume that the
convex angle determined by $u_1$ and $u_2$ does not contain any
other boundary segment with endpoint $X$. Furthermore, we may
assume that the segment $u_1$ is directed from $X$ to the other
endpoint. 

For $0<t<|u_1|$, let $A(t) \in u_1$ denote the point with $|A(t)-X|=t$ and 
let $A'(t)=P_4(A(t))$. There exists $\varepsilon >0$ such that for all 
$0<t<\varepsilon$ the points $A(t)$ are feasible, the points $A'(t)$ are feasible as 
well and lie on a common boundary segment. By our assumption, the convex 
angles between the ray $A(t)A'(t)$ and the segments $u_1, u_2$ directed from 
$X$ are both greater than $\frac{\pi}{2}$. It follows that if $t$ is 
sufficiently small, the tangent to the circle $\cC(A(t))$ at $A(t)$ 
intersects both segments $u_1, u_2$ and so does the circle $\cC(A(t))$, 
contradicting Claim~\ref{tvrz-sum}.
\end{proof}

An important consequence of Lemma~\ref{lem-uhel} is that no three boundary 
segments share a common endpoint. Hence, every connected component of the 
boundary is either an infinite piecewise linear curve, or a simple closed 
piecewise linear curve (i.e. the boundary of a simple polygon). We will 
call these cuves \emph{boundary components} or simply \emph{components}. 

\begin{defi}
Let $A$ be a point on the boundary. For $t\in\bbR$, let
$A(t)$ denote the point of the same boundary component as $A$,
such that the directed length of the part of the boundary starting
at $A$ and ending at $A(t)$ is equal to~$t$. $A(t)$ is clearly a
continuous function of $t$. If $A(t)$ is a feasible point, we let
$p_i(t)=P_i(A(t))$, for $i=0,\dotsc,5$.
\end{defi}

It is easy to see that the functions $p_i$ are continuous on a
sufficiently small neighborhood of every value of $t$ for which
$A(t)$ is a feasible point. Our next aim is to show that these
functions can be extended into continuous functions by suitably
defining the values of $p_i(t)$ when $A(t)$ is not feasible. It is
not obvious that the functions $p_i$ can be extended in this way:
the definition of $P_i(A(t))$ uses the Cartesian system whose
$x$-axis is parallel with the boundary segment containing $A(t)$.
Hence, if $A_1$ and $A_2$ are two feasible points belonging to two
distinct boundary segments of the same boundary component, it might
not be immediately clear that $P_i(A_1)$ belongs to the same
boundary component as $P_i(A_2)$. The next lemma shows that these
technical difficulties can be overcome.

\begin{lem}\label{lem-limita}
Let $A(t_0)$ be an infeasible point. For every $i=0,\dotsc,5$,
there is a point $P_i\in\cB$ such that

\[
\lim_{t\to t_0-}p_i(t)=P_i=\lim_{t\to t_0+} p_i(t)
\]
This means that if we define $p_i(t_0)=P_i$, then $p_i$ is continuous at $t_0$.
\end{lem}
\begin{proof}
It is sufficient to prove the lemma for $i=0$, because $p_i(t)$ is clearly a 
continuous function of $A(t)$ and $p_0(t)$. Since every boundary segment 
contains only finitely many infeasible points, we may choose a sufficiently 
small $\varepsilon>0$, such that for every $t$ from the open interval 
$(t_0-\varepsilon,t_0)$ the points $A(t)$ are feasible and they all belong 
to a single boundary segment $u_1$, and similarly, for every $t'\in 
(t_0,t_0+\varepsilon)$ the points $A(t')$ are feasible, and they belong to a 
single boundary segment $u_2$. If the segments $u_1$ and $u_2$ are distinct, 
then $A(t_0)$ is their common endpoint. Note that for $t 
\in(t_0-\varepsilon,t_0)$, the points $p_0(t)$ all belong to a single 
boundary segment $v_1$, otherwise some of the $A(t)$ would not be feasible. 
By Claim~\ref{tvrz-sum}, the segment $v_1$ is parallel and consistently 
oriented with $u_1$. Similarly, for $t'\in(t_0,t_0+\varepsilon)$ the points 
$p_0(t')$ belong to a single boundary segment $v_2$, parallel and consistently 
oriented with $u_2$. We do not know yet whether $v_1$ and $v_2$ appear 
consecutively on the same component of the boundary.

Let $B=\lim_{t\to t_0-} p_0(t)$ (clearly, the limit exists, because the 
points $\{p_0(t);\, t\!\in\!(t_0-\varepsilon,t_0)\}$ form an open segment whose 
endpoint is $B$). See Fig.~\ref{fig-limita}.

\begin{figure}[ht]
\begin{center}\includegraphics{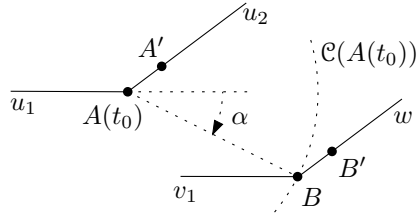}\end{center}
\caption{Illustration of the proof of 
Lemma~\ref{lem-limita}}\label{fig-limita}
\end{figure} 

For $t\in (t_0-\varepsilon,t_0)$, let us fix 
$\alpha\in(-\frac{\pi}{6},\frac{\pi}{6})$ such that $p_0(t)=P(\alpha,A(t))$, 
i.e., $\alpha$ is the (signed) measure of the angle between the segment $u_1$ 
and the segment $A(t)p_0(t)$. Note that $\alpha$ does not depend on the 
choice of $t$. The circle $\cC(A(t_0))$ intersects the boundary at $B$. Let 
$w$ be the boundary segment starting at $B$ and directed away from $B$. By 
Lemma~\ref{lem-uhel}, the convex angles determined by $v_1$ and $w$ and by 
$u_1$ and $u_2$ have size at least $\frac{2\pi}{3}$, which implies that the 
convex angle $\alpha'$ between $u_2$ and $BA(t_0)$ is acute and the convex 
angle between $w$ and $BA(t_0)$ is obtuse. Thus, for 
$t'\in(t_0,t_0+\varepsilon)$ the circle $\cC(A')$ (where $A'=A(t')$) 
intersects the segment $w$ at a point $B'=p_i(t')$. From Claim~\ref{tvrz-sum} 
it follows that $w$ is parallel to $u_2$. Also, the segment $A'B'$ is 
parallel to the segment $A(t_0)B$, which is in turn parallel to any of the 
segments $A(t)p_0(t)$, for $t\in(t_0-\varepsilon,t_0)$.

To finish the proof of this lemma, we need to show that $B'=p_0(t')$ (as 
opposed to $B'=p_i(t')$ for some $i\neq 0$), i.e., we need to prove that the 
angle $\alpha'$ determined by the segment $u_2$ and the segment $A'B'$ falls 
into the range $(-\frac{\pi}{6},\frac{\pi}{6})$. We have observed that 
$\alpha'\in(-\frac{\pi}{2},\frac{\pi}{2})$. This leaves us with the following 
three possibilities: either $B'=p_5(t')$, or $B'=p_1(t')$, or $B'=p_0(t')$. 
However, the former two possibilities are ruled out by the fact that the 
segment $w$ is oriented consistently with the segment $u_2$.
This concludes the proof.
\end{proof}

\begin{lem}\label{lem-posun}
Let $i\in\{0,\dotsc,5\}$, let $A\in\cB$ be an arbitrary boundary point.
All the unit segments of the form $A(t)p_i(t)$ have
the same slope, independently of the choice of $t$.
\end{lem}
\begin{proof}
The slope of $A(t)p_i(t)$ (as a function of $t$) is constant in a neighborhood
of every $t$ for which $A(t)$ is feasible. Moreover, this slope is a continuous
function of $t$, which follows from Lemma~\ref{lem-limita}. Hence the function
is constant on the whole range.
\end{proof}

Lemma~\ref{lem-posun} shows that every translation that maps a
feasible point $A$ to the point $P_i(A)$ also maps the boundary
component containing $A$ onto the boundary component containing
$P_i(A)$ (which may be the same component). Composing such
translations (or their inverses) we conclude that the translations
that send $P_i(A)$ to $P_j(A)$ have the same component-preserving
property.

For the proof of Lemma~\ref{lem-usek}, we will need a slight extension of 
Claim~\ref{tvrz-sum} to infeasible points:

\begin{cla}\label{cla-infeasible_body}
Let $A\in\cB$ be an arbitrary infeasible point.
\begin{enumerate}[(i)]
  
\item  
At each of the six points $P_0(A), P_1(A), \dots, P_5(A)$ the circle $\cC(A)$ 
properly crosses the corresponding boundary component, i.e., in a sufficiently 
small neighborhood of such point, the circle $\cC(A)$ separates the boundary 
component into two portions, one lying inside $\cC(A)$ and the other one 
lying outside $\cC(A)$. 
\item  
There are no more proper crossings of $\cC(A)$ with boundary components. 
(However, $\cC(A)$ may touch the boundary at some other points.)
\item 
The boundary components containing the points $P_0(A)$ and $P_3(A)$ have the 
same orientation as the component containing $A$, whereas the boundary 
components containing $P_1(A)$, $P_2(A)$, $P_4(A)$ and $P_5(A)$ have opposite 
orientation.
\end{enumerate}
\end{cla}

\begin{proof}
The first two statements follow from the fact that $\cC(A)$ has the same 
number of proper crossings with the boundary as the circle $\cC(A(t))$, where 
$A(t)$ is a feasible point sufficiently close to $A$. The third statement 
follows from Claim~\ref{tvrz-sum} applied to the point $A(t)$.
\end{proof}

\begin{lem}\label{lem-usek}
Let $A\in\cB$ be an arbitrary boundary point. For the sake of brevity, let us 
write $P_i$ instead of $P_i(A)$, $\cC$ instead of $\cC(A)$ and $\cD$ instead 
of $\cD(A)$ in the statement and proof of this lemma. The point $P_1$ belongs 
to the same boundary component as $P_2$, the point $P_0$ belongs to the same 
boundary component as $A$ and $P_3$, and the point $P_4$ belongs to the same 
boundary component as $P_5$. The four portions of the boundary that connect 
$P_1$ with $P_2$, $P_0$ with $A$, $A$ with $P_3$, and $P_4$ with $P_5$ are 
all translated copies of a single piecewise linear curve. These four portions 
of the boundary are all contained in the closed unit disc with center~$A$.  
\end{lem}
\begin{proof} It suffices to show that the boundary component that enters inside
$\cD$ at $P_1$ leaves $\cD$ at $P_2$. The rest of the statement
follows from Lemma~\ref{lem-posun}.

Let $\cL$ be the boundary component that contains $P_1$. Let us
follow $\cL$ from $P_1$ in the direction of its orientation, i.e.,
into the interior of the unit disc $\cD$, and let $X$ be the first
point where $\cL$ leaves $\cC$. We observe the following:
\begin{itemize}
\item
$X$ is neither $P_3$ nor $P_5$, because in these points, the
boundary is oriented into the interior of the disc $\cD$.
\item
$X$ is not the point $P_0$: if $X=P_0$, then the translation
$P_0\mapsto A$ would map the fragment of the boundary between
$P_1$ and $P_0$ onto a fragment directed from $P_2$ to $A$.
Similarly, the translation $P_1\mapsto A$ would map the fragment
$P_1P_0$ onto a fragment directed from $P_5$ to $A$. This is
impossible, because two different boundary fragments of equal
length cannot both end at $A$.
\item
$X$ is not $P_4$: if $X$ were equal to $P_4$, we would consider the boundary
component that enters into the interior of $\cC$ at the point $P_3$. Since this
boundary component cannot intersect the boundary fragment between $P_1$ and
$P_4$, it must leave the interior of $\cC$ at the point $P_2$. However, this is
symmetric to the previous case and leads to contradiction in the same way.
\item
Having excluded all other possibilities, we know that $X=P_2$.
\end{itemize}
Let $U$ denote the fragment of $\cL$ between $P_1$ and $P_2$. By definition, 
this fragment properly crosses $\cC$ only at its endpoints. Applying a 
symmetric argument, we find that the boundary fragment from $P_5$ to $P_4$ 
(which is a translated copy of $U$) properly crosses $\cC$ only in its 
endpoints. Translating $U$ appropriately, we obtain the boundary fragments 
connecting $P_3$ with $A$ and $A$ with $P_0$. This concludes the proof.
\end{proof}

From the previous lemmas, we readily obtain the following claim.
\begin{cla}\label{cla-nutne}
The condition (C3) of Theorem~\ref{thm-poly} implies the condition (C1).
\end{cla}                                                          
\begin{proof}                         
We check that the coloring $\chi$ satisfies the conditions of 
Definition~\ref{def-strip}. Let $\vec x$ denote the unit vector 
$\overrightarrow{AP_0}$ and let $\vec y$ be a unit vector orthogonal to $\vec 
x$. By Lemma~\ref{lem-usek}, every component of the boundary is a piecewise 
linear $\vec x$-periodic curve and if $\cL$ is a boundary component, then any 
other component is a translate of $\cL$ by an integral multiple of the vector 
$\overrightarrow{AP_1}=\frac{1}{2}\vec x+\frac{\sqrt{3}}{2}\vec y$. Let  
$\vec z$ denote this last vector and let $\cL_i=\cL_0+i\vec z$, $i\in\bbZ$, 
where $\cL_0$ is a boundary component chosen arbitrarily. We have 
$\cB=\bigcup_{i\in\bbZ}\cL_i$. The condition $(d)$ of 
Definition~\ref{def-strip} follows from Lemma~\ref{lem-usek}.       
\end{proof}

It remains to show that the condition (C1) implies (C2). This is the
easier part of the proof. In fact, we prove a more general claim:

\begin{thm}\label{thm-obarv}
Every zebra-like coloring has a twin that avoids the unit triangle.
\end{thm}                                                                
\begin{proof}                                                            
Let $\chi$ be a zebra-like coloring, let $\cL_i$, $\vec x$ and $\vec y$ 
be as in Definition~\ref{def-strip}. Let $\vec z=\frac{1}{2}\vec x+\frac{\sqrt{3}}{2}\vec 
y$. Let $\chi'$ be the twin coloring of 
$\chi$ such that the points of $\cL_i$ are black in $\chi'$ if $i$ is even and 
white if $i$ is odd.

Observe that by the definition of the coloring, the color of a point $P$ is 
equal to the color of $P+\vec x$ and different from the color of $P+\vec z$. 
Now assume that $ABC$ is a monochromatic unit triangle, wlog the three points 
are black. By the previous observation, no edge of the triangle forms an 
angle of size $\frac{\pi}{3}$ (or $\frac{2\pi}{3}$) with the vector $\vec x$. 
It follows that exactly one of the three edges (wlog the edge $AB$) forms 
with $\vec x$ an angle whose size falls into the range 
$(\frac{\pi}{3},\frac{2\pi}{3})$. 

We claim that the three points $A,B,C$ all belong to a single connected 
component of the black color: otherwise one of the two edges $AC$ and $BC$ 
would have to intersect (at least) two curves $\cL_i$ and $\cL_{i+1}$. By the 
definition of the coloring, the distance between the two points of intersection 
is greater than~1, contradicting the fact that $ABC$ is a unit triangle.

We now deduce that $\|AB\|<1$: let $\ell$ be the line containing the segment 
$AB$. Note that the line $\ell$, as well as any other line not parallel with 
$\vec x$, must intersect all the curves $\cL_i$. Let $A'B'$ be the segment 
obtained as the convex hull of the intersection of $\ell$ with the closure of 
the black component containing $A$ and $B$. By the definition of the coloring, 
$\|A'B'\|\le 1$. Moreover, since the two points $A'$ and $B'$ belong to two 
adjacent boundary curves $\cL_i$ and $\cL_{i+1}$, they have different colors. 
Hence, the segment $AB$ is a proper subset of the segment $A'B'$, and  
$\|AB\|<1$. This shows that $ABC$ is not a unit triangle\thinspace ---\thinspace a~contradiction.
\end{proof}

This concludes the proof of Theorem~\ref{thm-poly}. Next, we
present a simple corollary, which shows that every polygonal coloring of
the plane contains any nonequilateral triangle.

\subsection{Nonequilateral triangles}
The following result is a direct consequence of Theorem~\ref{thm-poly}, by an 
easy modification of the proof of Lemma~\ref{lem-osm}. 

\begin{thm}\label{thm-osm}
Let $XYZ$ be a nonequilateral triangle, let $\chi$ be a polygonal coloring. 
There is a monochromatic copy $X'Y'Z'$ of the configuration $XYZ$, such that 
none of the three points $X',Y'$ and $Z'$ belongs to the boundary of $\chi$.
\end{thm}
\begin{proof}
Let $a, b$ and $c$ be the lengths of the three edges of $XYZ$. Wlog, assume 
that $a\neq b$. From Theorem~\ref{thm-poly} it follows that no polygonal 
coloring can simultaneously avoid copies of equilateral triangles of two 
different sizes. Hence, we may assume that $\chi$ contains a monochromatic 
equilateral triangle $ABC$ with edges of length $a$ whose vertices avoid the 
boundary of $\chi$. Assume that the three points $A$, $B$ and $C$ are all 
black. Consider the configuration of eight points on Fig.~\ref{fig-osm}. As 
discussed in the proof of the first part of Lemma~\ref{lem-osm}, every 
coloring of the five points $D$, $A'$, $B'$, $C'$ and $D'$ yields a 
monochromatic $(a,b,c)$-triangle. Furthermore, we may assume that the eight 
points all avoid the boundary of $\chi$, otherwise we might shift the 
configuration slightly to move the points away from the boundary, without changing the 
color of $ABC$ (recall that $A, B$ and $C$ already belong to the interior of 
the black color). This concludes the proof.
\end{proof}  

\section{Concluding remarks}

The Conjecture~\ref{con-slaba} remains wide open, despite the indirect 
support from the results of this paper, as well as from earlier research. It 
might happen that the validity of this conjecture would depend on the 
particular choice of set-theoretical axioms. Such issues do not arise in this 
paper, since our proof techniques are very elementary. Unfortunately, these 
elementary techniques do not offer much hope for broad generalizations. It 
might nevertheless be possible to extend our results about polygonal colorings 
to some broader class of colorings, e.g., the colorings by monochromatic 
regions bounded by continuous curves. Colorings of this kind have already been 
studied in the context of the related problem of the chromatic number of the 
plane (see \cite{woo}).

The zebra-like colorings provide a hitherto unknown example of 
colorings that avoid an equilateral triangle. We are not aware of any other 
examples of colorings avoiding a given triangle, but we do not dare to make any 
conjectures about the uniqueness of our construction, because our understanding 
of non-polygonal colorings is rather limited.        

\section*{Acknowledgments}
We appreciate the useful discussions with Zden\v ek Dvo\v r\'ak, Jan 
Kratochv\'\i l, Martin Tancer, Pavel Valtr and Tom\'a\v s Vysko\v cil.

\end{document}